# COUPLING WITH THE STATIONARY DISTRIBUTION AND IMPROVED SAMPLING FOR COLORINGS AND INDEPENDENT SETS

By Thomas P. Hayes[1] and Eric Vigoda[2]

*University of California at Berkeley and Georgia Institute of Technology*

We present an improved coupling technique for analyzing the mixing time of Markov chains. Using our technique, we simplify and extend previous results for sampling colorings and independent sets. Our approach uses properties of the stationary distribution to avoid worst-case configurations which arise in the traditional approach.

As an application, we show that for $k/\Delta > 1.764$, the Glauber dynamics on $k$-colorings of a graph on $n$ vertices with maximum degree $\Delta$ converges in $O(n \log n)$ steps, assuming $\Delta = \Omega(\log n)$ and that the graph is triangle-free. Previously, girth $\geq 5$ was needed.

As a second application, we give a polynomial-time algorithm for sampling weighted independent sets from the Gibbs distribution of the hard-core lattice gas model at fugacity $\lambda < (1-\varepsilon)e/\Delta$, on a regular graph $G$ on $n$ vertices of degree $\Delta = \Omega(\log n)$ and girth $\geq 6$. The best known algorithm for general graphs currently assumes $\lambda < 2/(\Delta - 2)$.

**1. Introduction.** The coupling method is an elementary yet powerful technique for bounding the rate of convergence of a Markov chain to its stationary distribution. Traditionally, the coupling technique has been a standard tool in probability theory (e.g., [1, 4, 14]) and statistical physics (e.g., [15]). More recently, it has yielded significant results in theoretical computer science [3, 5, 9, 10, 11, 16, 17]. We refine the coupling method, and as a consequence, improve and simplify recent results on randomly sampling colorings and weighted independent sets.

Consider a Markov chain on a finite state space $\Omega$, that has a unique stationary distribution $\pi$. A (one-step) *coupling* specifies, for every pair of

Received March 2004; revised November 2005.

[1]Supported in part by an NSF Postdoctoral Fellowship.
[2]Supported in part by NSF Grant CCF-04-55666.

states $(X_t, Y_t) \in \Omega^2$, a distribution for $(X_{t+1}, Y_{t+1})$ such that $X_t \to X_{t+1}$, and similarly $Y_t \to Y_{t+1}$, behave according to the Markov chain.

Let $\rho$ denote an arbitrary integer-valued metric on $\Omega$, where $\operatorname{diam}(\Omega)$ denotes the length of the longest path. For $\varepsilon > 0$, we say a pair $(x, y) \in \Omega^2$ is $\varepsilon$ *distance-decreasing* if there exists a coupling such that

$$\mathbf{E}(\rho(X_1, Y_1) | X_0 = x, Y_0 = y) < (1 - \varepsilon)\rho(x, y).$$

The coupling theorem says that if *every pair* $(x, y)$ is distance-decreasing, then the Markov chain mixes rapidly:

THEOREM 1.1 (cf. [1]). *Let $\varepsilon > 0$ and suppose every $(x, y) \in \Omega^2$ is $\varepsilon$ distance-decreasing. Let $X_0 \in \Omega$, $\delta > 0$ and*

$$T \geq \frac{\ln(\operatorname{diam}(\Omega)/\delta)}{\varepsilon}.$$

*Then $\|X_T - \pi\|_{\mathrm{TV}} \leq \delta$.*

For a pair of distributions $\mu, \nu$ on a space $\Omega$, the variation distance metric is defined as

$$\|\mu - \nu\|_{\mathrm{TV}} = \tfrac{1}{2} \sum_{x \in \Omega} |\mu(x) - \nu(x)|.$$

We will often abuse notation, writing $\|X - \nu\|_{\mathrm{TV}}$ for the variation distance between the distribution of a random variable $X$ and a distribution $\nu$.

Our first *coupling with stationarity* theorem does not require every pair of states to be distance-decreasing. Instead, we only require that *most states $x$* be distance-decreasing with *every $y$*.

THEOREM 1.2. *Let $\varepsilon > 0$. Suppose $S \subseteq \Omega$ such that every $(x, y) \in S \times \Omega$ is $\varepsilon$ distance-decreasing, and*

$$\pi(S) \geq 1 - \frac{\varepsilon}{16 \operatorname{diam}(\Omega)}.$$

*Let $X_0 \in \Omega$, $\delta > 0$ and*

$$T \geq \frac{\lceil \ln(32 \operatorname{diam}(\Omega)) \rceil \lceil \ln(1/\delta) \rceil}{\varepsilon}.$$

*Then $\|X_T - \pi\|_{\mathrm{TV}} \leq \delta$.*

We will apply Theorem 1.2 to improve results on randomly sampling colorings; see Section 1.1. Theorem 1.2 is proved in Section 3.

If we only require that *most pairs* of states be distance-decreasing, we can prove rapid mixing under the additional assumption that the initial distribution is a "warm start," as defined by Kannan, Lovasz and Simonovitz [13] for random walks in convex bodies. A distribution $X_0$ on $\Omega$ is said to be a warm start (with respect to $\pi$) if

$$\text{for all } z \in \Omega \qquad \mathbf{Pr}(X_0 = z) \leq 2\pi(z).$$

THEOREM 1.3. *Let $\varepsilon > 0$. Suppose $S \subseteq \Omega$ such that every $(x, y) \in S^2$ is $\varepsilon$ distance-decreasing. Let*

$$\pi(S) > 1 - \frac{\varepsilon \delta}{6 \operatorname{diam}(\Omega)} \quad and \quad T \geq \frac{\ln(2 \operatorname{diam}(\Omega)/\delta)}{\varepsilon}.$$

*Then if $X_0$ is a warm start to $\pi$,*

$$\|X_T + \pi\|_{\mathrm{TV}} < \delta.$$

We will apply Theorem 1.3 to improve results on randomly sampling independent sets (see Section 1.2). Theorem 1.3 is proved in Section 3.

1.1. *Randomly sampling colorings.* For a graph $G = (V, E)$ with maximum degree $\Delta$, the Glauber dynamics (heat-bath) is a simple Markov chain whose stationary distribution is uniformly distributed over proper $k$-colorings of $G$. Let $\Omega = [k]^V$, where $[k] = \{1, 2, \ldots, k\}$. From $X_t \in \Omega$, the evolution $X_t \to X_{t+1}$ is defined as follows:

- Choose $v$ uniformly at random from $V$.
- For all $w \neq v$, set $X_{t+1}(w) = X_t(w)$.
- Choose $X_{t+1}(v)$ uniformly at random from $[k] \setminus X_{t+1}(N(v))$, where $N(v)$ is the set of neighbors of $v$. In words, the new color for $v$ is randomly chosen from those colors not appearing in the neighborhood of $v$.

The latest results on randomly sampling colorings, beginning with [5], use the following "burn-in" method for the analysis of the Glauber dynamics. After the Markov chain evolves for a sufficient number of steps, the so-called burn-in period, the coloring has certain "local uniformity" properties with high probability. Moreover, these properties persist for a polynomial number of steps. Consequently, to prove rapid mixing it suffices to prove there is a distance-decreasing coupling for every pair of states satisfying the local uniformity properties. Earlier works, for example, [11], analyze the worst-case pair of states and hence rely upon Theorem 1.1. Using the burn-in approach led to many significant improvements [5, 9, 10, 16] since it avoids the worst-case pair of states in the coupling analysis.

Proving that the local uniformity properties appear for the Glauber dynamics is very difficult. Roughly speaking, vertex colors are not independent, and their correlation has to be bounded. Dyer and Frieze [5] and Molloy [16] used the method of "paths of disagreement." Hayes [9] used a more sophisticated method of "conditional independence." As a by-product of these results, it follows immediately that a uniformly random coloring has the local uniformity property with high probability.

Directly proving that a uniformly random coloring has these local uniformity properties is much easier than for colorings generated by the Glauber dynamics. Our upcoming Theorem 1.4 highlights this simplicity. By "coupling with stationarity" (Theorem 1.2), we are able to improve the main result of [9] with considerably less work than the original.

THEOREM 1.4. *Let $\alpha = 1.763\ldots$ denote the solution to $x = \exp(1/x)$. Let $0 < \zeta < 1$. Let $G$ be a triangle-free graph on n vertices having maximum degree $\Delta$, let $k \geq \max\{(1+\zeta)\alpha\Delta, 288\ln(96n^3/\zeta)/\zeta^2\}$ and let $X_0$ be any k-coloring of G. Then for every $\delta > 0$, after $T \geq 6n\lceil\ln(32n)\rceil\lceil\ln(1/\delta)\rceil/\zeta$ steps of the Glauber dynamics,*

$$\|X_T - \pi\|_{\mathrm{TV}} \leq \delta.$$

Earlier versions of Theorem 1.4 appeared in [5] and [9]. Both needed higher girth [$\Omega(\log \Delta)$ and $\geq 5$, resp.] and had considerably more difficult proofs. The proof of Theorem 1.4 is presented in Section 2.

Hayes and Vigoda [10] have proved $O(n \log n)$ mixing time for $k > (1+\varepsilon)\Delta$ for all $\varepsilon > 0$, assuming $\Delta = \Omega(\log n)$ and girth $\geq 9$. Recently, Dyer, Frieze, Hayes and Vigoda [6] reduced the condition on $\Delta$ to a sufficiently large constant, assuming $k/\Delta > 1.489\ldots$ and girth $\geq 6$.

Our results for graph colorings are syntactically similar to recent work by Goldberg, Martin and Paterson [8], who also examine $k$-colorings of triangle-free graphs for $k/\Delta > 1.763\ldots$. Their focus is on proving, for random colorings, that correlation between the colors assigned to two vertices decays exponentially fast with the distance between the pair. More precisely, they are proving a variant of a so-called strong spatial mixing property holds. For amenable graphs, strong spatial mixing is closely related to rapid mixing of the Glauber dynamics (cf. [7]).

They prove their version of strong spatial mixing holds for every triangle-free graph for $k/\Delta > 1.763\ldots$ for all $\Delta$. Their proof and our proof utilize similar local uniformity properties of the stationary distribution. However, in their setting it suffices for the properties to hold in expectation, whereas in the analysis of the dynamics it appears essential for the properties to hold with high probability. Contrasting our results and proofs with theirs highlights the differences between strong spatial mixing and rapid mixing of the Glauber dynamics for general graphs.

1.2. *Randomly sampling independent sets.* Given a graph $G = (V, E)$ and a fugacity $\lambda > 0$, the hard-core lattice gas model (see [2]) is defined on the set

$\Omega$ of independent sets of $G$. The weight of $X \subseteq V$, where $X \in \Omega$, is defined as
$$w(X) = \lambda^{|X|}.$$
We are interested in sampling from the (Boltzmann) Gibbs distribution $\pi$ on $\Omega$ where $\pi(X) = w(X)/Z$ and
$$Z = Z(G, \lambda) = \sum_{X \in \Omega} w(X)$$
is the partition function.

As with colorings, the Glauber dynamics for the hard-core model updates a random vertex at each step. From $X_t \in \Omega$, the transition $X_t \to X_{t+1}$ is defined by:

- Choose a vertex $v$ uniformly at random from $V$.
- Set
$$X' = \begin{cases} X_t \cup v, & \text{with probability } \lambda/(1+\lambda), \\ X_t \setminus v, & \text{with probability } 1/(1+\lambda). \end{cases}$$
- If $X' \in \Omega$, set $X_{t+1} = X'$; otherwise set $X_{t+1} = X_t$.

It is clear that the Glauber dynamics is reversible and ergodic, and the unique stationary distribution is $\pi$.

The latest result for general graphs is $O(n \log n)$ mixing time for $\lambda < 2/(\Delta - 2)$ by Vigoda [18]. It is widely believed the chain mixes rapidly for all
$$\lambda < \frac{(\Delta - 1)^{\Delta - 1}}{(\Delta - 2)^{\Delta}} \sim e/\Delta.$$
Applying Theorem 1.3, we prove there exists an efficient algorithm which reaches the above threshold for large-degree regular graphs with girth at least 6. To guarantee the warm-start condition, we use a simulated annealing algorithm similar to Jerrum, Sinclair and Vigoda's algorithm for estimating the permanent [12].

THEOREM 1.5. *For all $\varepsilon > 0$, there exists $C > 0$ such that for every $\delta > 0$ the following holds. Let $G$ be a $\Delta$-regular graph on $n$ vertices, where $girth(G) \geq 6$, and $\Delta \geq C \log(n/\delta)$. Let $\lambda \leq (1-\varepsilon)e/\Delta$. Then there is an algorithm which outputs a random independent set of $G$ within variation distance $\delta$ of the Gibbs distribution at fugacity $\lambda$, with running time polynomial in $n, 1/\varepsilon$ and $\log(1/\delta)$.*

Theorem 1.5 is proved in Section 4. The degree restriction prevents us from boosting the above sampling scheme to arbitrarily close distances to the stationary distribution, and also from applying the standard reduction from approximating the partition function to sampling from the Gibbs distribution.

## 2. Sampling colorings.

2.1. *Local uniformity property.* We will exploit a nice "local uniformity" property of the uniform distribution on proper $k$-colorings. This property was first used by Dyer and Frieze [5] to improve rapid mixing results of the Glauber dynamics.

For any triangle-free graph we will show an easy lower bound on the expected number of available colors for an arbitrary vertex $v$, where "available colors" refers to those colors not appearing in the neighborhood of $v$. When $k = \Omega(\log n)$, we can even prove that, with high probability, every vertex has essentially this many available colors.

For any coloring $X$ which satisfies our lower bound on available colors, it is easy to show that for *every* coloring $Y$, the pair $(X, Y)$ is distance-decreasing under the natural "greedy" one-step coupling of Jerrum [11]. This allows us to apply our first "coupling with stationarity" Theorem 1.2 to prove Theorem 1.4.

Throughout, we will use the notation

$$A(X, v) := [k] \setminus X(N(v))$$

to denote the set of available colors for a vertex $v$ under a $k$-coloring $X$ [here $N(v) := \{w \in V | w \sim v\}$ denotes the set of neighbors of vertex $v$].

LEMMA 2.1. *Let $G = (V, E)$ be a triangle-free graph with maximum degree $\Delta$, let $0 < \beta \leq 1$ and $k \geq \Delta + 2/\beta$. Let $X$ be a random $k$-coloring of $G$. Then*

$$\mathbf{Pr}(\exists v \in V, |A(X, v)| < k(e^{-\Delta/k} - \beta)) \leq n e^{-\beta^2 k/8}.$$

PROOF. Let $v \in V$. By definition,

$$|A(X, v)| = \sum_{j \in [k]} \prod_{w \in N(v)} (1 - X_{j,w}), \tag{1}$$

where $X_{j,w}$ is the indicator variable for the event $\{X(w) = j\}$.

Henceforth, we will condition on the values of $X$ on $V \setminus N(v)$; denote this conditional information by $\mathcal{F}$. Conditioned on $\mathcal{F}$, since $G$ is triangle-free, the random variables $X(w), w \in N(v)$, are fully independent and each $X(w)$ is uniform over the set $A(X, w)$. This allows us to write

$$\mathbf{E}(|A(X, v)||\mathcal{F}) = \sum_{j \in [k]} \prod_{w \in N(v)} (1 - \mathbf{E}(X_{j,w}|\mathcal{F}))$$
$$= \sum_{j \in [k]} \prod_{\substack{w \in N(v) \\ A(X,w) \ni j}} \left(1 - \frac{1}{|A(X, w)|}\right).$$

Applying the arithmetic–geometric mean inequality, this implies

$$\mathbf{E}(|A(X,v)||\mathcal{F}) \geq k \prod_{\substack{j\in[k] \ w\in N(v) \\ A(X,w)\ni j}} \left(1 - \frac{1}{|A(X,w)|}\right)^{1/k}$$

$$= k \prod_{w\in N(v)} \prod_{j\in A(X,w)} \left(1 - \frac{1}{|A(X,w)|}\right)^{1/k}$$

$$= k \prod_{w\in N(v)} \left(1 - \frac{1}{|A(X,w)|}\right)^{|A(X,w)|/k}$$

$$\geq k \prod_{w\in N(v)} \left(\frac{1 - 1/|A(X,w)|}{e}\right)^{1/k}$$

where the last step follows from the inequality $(1-p)^{1/p} \geq (1-p)/e$, which holds for all $0 < p \leq 1$. Since we are assuming $k \geq \Delta + 2/\beta$, it follows that

$$\mathbf{E}(|A(X,v)||\mathcal{F}) \geq ke^{-\Delta/k}\left(1 - \frac{\beta}{2}\right)^{\Delta/k}$$

$$\geq k\left(e^{-\Delta/k} - \frac{\beta}{2}\right).$$

Since the colors $X(w), w \in N(v)$, are fully independent, conditioned on $\mathcal{F}$, and since $|A(X,w)|$ is a Lipschitz function of these colors with Lipschitz constant 1, it follows by Chernoff's bounds that

$$\mathbf{Pr}(|A(X,v)| \leq k(e^{-\Delta/k} - \beta)) \leq e^{-\beta^2 k/8}.$$

Taking a union bound over $v \in V$ completes the proof. $\square$

REMARK 2.2. A stronger form of Lemma 2.1 was proved by Hayes [9], who replaced the assumption $k \geq \Delta + 2/\beta$ by $k \geq \Delta + 2$ with slightly worse constants in the error probability bound.

2.2. *Most colorings are distance-decreasing.* We now present a simple sufficient condition for a pair of colorings to be distance-decreasing. We use Jerrum's one-step coupling of the Glauber dynamics on $k$-colorings; see [11]. Each chain chooses the same vertex to recolor at every step. We then maximize the probability the chains choose the same new color for the updated vertex. Under this coupling, for the purposes of proving rapid mixing it suffices for one of the chains to have the local uniformity property considered in Lemma 2.1.

LEMMA 2.3. *Let $0 < \beta < 1$, and suppose $X \in \Omega$ satisfies, for every $v \in V(G)$,*

$$|A(X,v)| \geq \frac{\Delta}{1-\beta}.$$

*Then, for every $Y \in \Omega$, the pair $(X,Y)$ is $\beta/n$ distance-decreasing.*

PROOF. We need to prove, for every $Y \in \Omega$,

$$\mathbf{E}(\rho(X_1, Y_1) | X_0 = X, Y_0 = Y) \leq \left(1 - \frac{\beta}{n}\right) \rho(X, Y),$$

where in this case $\rho$ denotes Hamming distance on graph colorings. Let us fix a coloring $Y \in \Omega$, and condition throughout on the event $\{X_0 = X, Y_0 = Y\}$.

Let $v$ denote the vertex randomly selected for recoloring at time 1. Jerrum's coupling maximizes the probability the chains choose the same color for $v$. For a color $c$ available to $v$ in both chains, that is, $c \in A(X,v) \cap A(Y,v)$, we simultaneously color $v$ with $c$ in both $X$ and $Y$ with probability $\min\{1/|A(X,v)|, 1/|A(Y,v)|\}$. With the remaining probabilities the coupled color choices are arbitrary. Hence,

$$\mathbf{Pr}(X_1(v) = Y_1(v) = c | v) = \frac{1}{\max\{|A(X,v)|, |A(Y,v)|\}}.$$

We now bound the probability the chains recolor $v$ to a different color in the two chains,

$$\mathbf{Pr}(X_1(v) \neq Y_1(v) | v) = 1 - \frac{|A(X,v) \cap A(Y,v)|}{\max\{|A(X,v)|, |A(Y,v)|\}}$$

$$= \frac{\max\{|A(X,v)|, |A(Y,v)|\} - |A(X,v) \cap A(Y,v)|}{\max\{|A(X,v)|, |A(Y,v)|\}}.$$

Hence,

(2) $$\mathbf{Pr}(X_1(v) \neq Y_1(v) | v) \leq \frac{|\{u \in N(v) | X(u) \neq Y(u)\}|}{\max\{|A(X,v)|, |A(Y,v)|\}}.$$

Finally, we bound the expected distance after the transition,

$$\mathbf{E}(\rho(X_1, Y_1)) = \sum_{w \in V} \mathbf{Pr}(X_1(w) \neq Y_1(w))$$

$$= \sum_{w \in V} \mathbf{Pr}(v \neq w \text{ and } X(w) \neq Y(w))$$

$$+ \sum_{w \in V} \mathbf{Pr}(v = w \text{ and } X_1(w) \neq Y_1(w))$$

$$= \frac{n-1}{n}\rho(X,Y) + \sum_{w \in V} \mathbf{Pr}(v = w \text{ and } X_1(w) \neq Y_1(w))$$

$$\leq \frac{n-1}{n}\rho(X,Y) + \frac{1}{n}\sum_{w \in V} \frac{|\{u \in N(w)|X(u) \neq Y(u)\}|}{\max\{|A(X,w)|, |A(Y,w)|\}}$$

$$\leq \frac{n-1}{n}\rho(X,Y) + \frac{1}{n}\frac{\rho(X,Y)\Delta}{\Delta/(1-\beta)}$$

$$= \left(1 - \frac{\beta}{n}\right)\rho(X,Y).$$

The first inequality above holds by (2), and the second inequality uses the fact that $\Delta$ is the maximum degree and our assumption that $|A(X,w)| \geq \Delta/(1-\beta)$ for all $w \in V$. □

We now present the proof of Theorem 1.4.

PROOF OF THEOREM 1.4. Define $\beta = \zeta/6$. Since by hypothesis, $k \geq (1+\zeta)\alpha\Delta > 3\Delta/2$ and $k \geq 288\ln(96n^3/\zeta)/\zeta^2 = 8\ln(16n^3/\beta)/\beta^2 > 8/\beta$, it follows that either $\Delta \geq 4/\beta$, and so $k \geq 3\Delta/2 \geq \Delta + 2/\beta$, or $\Delta < 4/\beta$, and so $k \geq 8/\beta \geq \Delta + 2/\beta$. Thus in either case $k \geq \Delta + 2/\beta$. Let

$$S = \{X \in \Omega : (\forall v \in V)[|A(X,v)| \geq k(e^{-\Delta/k} - \beta)]\}.$$

Lemma 2.1 together with the hypothesis $k \geq 8\ln(16n^3/\beta)/\beta^2$ and the fact $\text{diam}(\Omega) = n$, imply

$$\pi(S) \geq 1 - ne^{-\beta^2 k/8}$$

$$\geq 1 - \frac{\beta}{16n^2}$$

$$= 1 - \frac{\varepsilon}{16\,\text{diam}(\Omega)},$$

where $\varepsilon = \beta/n$.

Recalling that $\alpha = 1.763\ldots < 2$ and that $0 < \zeta = 6\beta \leq 1$, it can be verified by elementary algebra that

$$(1+\zeta)(1-\beta)(1-\alpha\beta) \geq 1,$$

and hence that

$$\frac{\Delta}{1-\beta} \leq \frac{k}{\alpha(1+\zeta)(1-\beta)}$$

$$\leq k\left(\frac{1}{\alpha} - \beta\right)$$

$$= k(\exp(-1/\alpha) - \beta) \qquad \text{by definition of } \alpha$$

$$< k(\exp(-\Delta/k) - \beta).$$

It follows that, by Lemma 2.3, every pair $(X,Y) \in S \times \Omega$ is $\varepsilon$ distance-decreasing, where $\varepsilon = \beta/n$. Applying Theorem 1.2 yields the desired result. □

**3. Coupling with stationarity.** In this section we will prove our coupling with stationarity theorems (Theorems 1.2 and 1.3). These will follow as corollaries of the following more general theorem about couplings which "usually" decrease distances. For an event $\mathcal{G}$, let $\mathbf{1}_{\mathcal{G}}$ denote the indicator variable whose value is 1 if $\mathcal{G}$ occurs and 0 otherwise.

THEOREM 3.1. *Let $X_0, \ldots, X_T, Y_0, \ldots, Y_T$ be coupled Markov chains such that, for every $0 \leq t \leq T-1$,*

$$\mathbf{Pr}((X_t, Y_t) \text{ is not } \varepsilon \text{ distance-decreasing}) \leq \delta.$$

*Then*

$$\mathbf{Pr}(X_T \neq Y_T) \leq ((1-\varepsilon)^T + \delta/\varepsilon) \operatorname{diam}(\Omega).$$

PROOF. For $0 \leq t \leq T-1$, let $\mathcal{G}(t)$ denote the event that $(X_t, Y_t)$ is $\varepsilon$ distance-decreasing. Observe that

$$\begin{aligned}
&\mathbf{E}(\rho(X_{t+1}, Y_{t+1}) - (1-\varepsilon)\rho(X_t, Y_t)) \\
&= \mathbf{E}((\rho(X_{t+1}, Y_{t+1}) - (1-\varepsilon)\rho(X_t, Y_t))\mathbf{1}_{\mathcal{G}(t)}) \\
&\quad + \mathbf{E}((\rho(X_{t+1}, Y_{t+1}) - (1-\varepsilon)\rho(X_t, Y_t))\mathbf{1}_{\overline{\mathcal{G}(t)}}) \\
&\leq \mathbf{E}((\rho(X_{t+1}, Y_{t+1}) - (1-\varepsilon)\rho(X_t, Y_t))\mathbf{1}_{\overline{\mathcal{G}(t)}}) \\
&\leq \delta \mathbf{E}(\rho(X_{t+1}, Y_{t+1})) \\
&\leq \delta \operatorname{diam}(\Omega).
\end{aligned}$$

This can be rewritten as

$$(3) \qquad \mathbf{E}(\rho(X_{t+1}, Y_{t+1})) \leq (1-\varepsilon)\mathbf{E}(\rho(X_t, Y_t)) + \delta \operatorname{diam}(\Omega).$$

Hence,

$$\begin{aligned}
\mathbf{E}(\rho(X_1, Y_1)) &\leq (1-\varepsilon)\mathbf{E}(\rho(X_0, Y_0)) + \delta \operatorname{diam}(\Omega) \\
&\leq (1-\varepsilon+\delta) \operatorname{diam}(\Omega).
\end{aligned}$$

And, for all $t \geq 0$, from (3), it follows by induction that

$$\begin{aligned}
\mathbf{E}(\rho(X_t, Y_t)) &\leq ((1-\varepsilon)^t + \delta(1-(1-\varepsilon)^t)/\varepsilon) \operatorname{diam}(\Omega) \\
&< ((1-\varepsilon)^t + \delta/\varepsilon) \operatorname{diam}(\Omega).
\end{aligned}$$

Since $\rho$ is integer-valued, we obtain the desired conclusion by Markov's inequality:

$$\begin{aligned}\mathbf{Pr}(X_T \neq Y_T) &= \mathbf{Pr}(\rho(X_T, Y_T) \geq 1) \\ &\leq \mathbf{E}(\rho(X_T, Y_T)) \\ &< ((1-\varepsilon)^T + \delta/\varepsilon) \operatorname{diam}(\Omega). \end{aligned} \qquad \square$$

We now present the proofs of our coupling with stationarity theorems.

PROOF OF THEOREM 1.2. Let $X_0$ be arbitrary, and let $Y_0$ be distributed according to $\pi$. Generate $X_1, \ldots, X_T, Y_1, \ldots, Y_T$ using the given coupling with initial states $X_0, Y_0$. Note that, for every $t \geq 0$, $Y_t$ is distributed according to $\pi$, and so, since every element of $\Omega \times S$ is $\varepsilon$ distance-decreasing,

$$\begin{aligned}\mathbf{Pr}((X_t, Y_t) \text{ is not } \varepsilon \text{ distance-decreasing}) &\leq \mathbf{Pr}(Y_t \notin S) \\ &= 1 - \pi(S).\end{aligned}$$

Let $T' = \lceil \ln(32 \operatorname{diam}(\Omega))/\varepsilon \rceil$. We first show that after $T'$ steps, we are within distance $1/8$ of stationarity. Then the standard boosting argument implies that after $T = T' \lceil \ln(1/\delta) \rceil$ steps, we are within distance $\delta$ of stationarity.

Applying Theorem 3.1, and noting again that $Y_{T'} \sim \pi$, we have

$$\begin{aligned}\|X_{T'} - \pi\| &\leq \mathbf{Pr}(X_{T'} \neq Y_{T'}) \\ &\leq ((1-\varepsilon)^{T'} + (1-\pi(S))/\varepsilon) \operatorname{diam}(\Omega) \\ &\leq \left(\exp(-\varepsilon T') + \frac{1}{16 \operatorname{diam}(\Omega)}\right) \operatorname{diam}(\Omega) \\ &\leq 1/8.\end{aligned}$$

Since the above holds for all $X_0 \in \Omega$, the triangle inequality implies that for all pairs of states $W_0, Z_0 \in \Omega^2$,

$$\|W_{T'} - Z_{T'}\|_{\mathrm{TV}} \leq 1/4 < 1/e.$$

Since there always exists a $T'$-step coupling which achieves the variation distance, the above can be boosted as follows (see [1]). We consider the $T$-step coupling generated by concatenating this $T'$-step coupling $\lceil \ln(1/\delta) \rceil$ times. Then

$$\begin{aligned}\|W_T - \pi\|_{\mathrm{TV}} &\leq \max_{Z_0} \mathbf{Pr}(W_T \neq Z_T) \\ &\leq \prod_{i=1}^{\lceil \ln(1/\delta) \rceil} \mathbf{Pr}(W_{iT'} \neq Z_{iT'} | W_{(i-1)T'} \neq Z_{(i-1)T'}) \\ &\leq (1/e)^{\lceil \ln(1/\delta) \rceil} \leq \delta.\end{aligned}$$

□

PROOF OF THEOREM 1.3. Let $X_0$ be a warm start to $\pi$, and let $Y_0 \sim \pi$. It follows from the definition of stationarity that $Y_t$ is distributed according to $\pi$ for all $t > 0$. We observe further that $X_t$ is a warm start for all $t > 0$ since, assuming $X_{t-1}$ is a warm start, we have, for every $x \in \Omega$,

$$\mathbf{Pr}(X_t = x) = \sum_{x' \in \Omega} \mathbf{Pr}(X_{t-1} = x') P(x', x)$$
$$\leq \sum_{x' \in \Omega} 2\pi(x') P(x', x)$$
$$= 2\pi(x).$$

Since every element of $S \times S$ is $\varepsilon$ distance-decreasing, it follows that

$$\mathbf{Pr}((X_t, Y_t) \text{ is not } \varepsilon \text{ distance-decreasing}) \leq \mathbf{Pr}(X_t \notin S) + \mathbf{Pr}(Y_t \notin S)$$
$$\leq 3(1 - \pi(S)).$$

Applying Theorem 3.1, and noting again that $Y_T \sim \pi$, we have

$$\|X_T - \pi\|_{\text{TV}} \leq \mathbf{Pr}(X_T \neq Y_T)$$
$$\leq \left((1-\varepsilon)^T + \frac{3(1-\pi(S))}{\varepsilon}\right) \text{diam}(\Omega)$$
$$< \delta,$$

for $T > \ln(2 \operatorname{diam}(\Omega)/\delta)/\varepsilon$. □

**4. Independent sets.** In this section we will prove Theorem 1.5, presenting an algorithm based on simulated annealing, which allows sampling from the Gibbs distribution for the hard-core model at fugacity $\lambda$, approaching the (believed) critical threshold $\lambda = e/\Delta$. Our algorithm, which we will present shortly, relies on the efficient convergence of the Glauber dynamics, given a warm start. More formally, we will require the following result, which is analogous to Theorem 1.4 for graph colorings.

LEMMA 4.1. *Let $\zeta, \delta > 0$. Let $G$ be a $\Delta$-regular graph on $n$ vertices having girth at least 6, where $\Delta \geq 320000 \ln(144 n^3/\zeta\delta)/\zeta^4$, and let $\lambda \leq (1-\zeta)e/\Delta$ and $T \geq 8n \ln(2n/\delta)/\varepsilon$. Let $X_0$ be a warm start to $\pi$, the Gibbs distribution for the hard-core lattice gas model at fugacity $\lambda$. Then after $T$ steps of the Glauber dynamics,*

$$\|X_T - \pi\|_{\text{TV}} \leq \delta.$$

We first present the algorithm of Theorem 1.5, which utilizes Lemma 4.1. We then prove Lemma 4.1.

*Simulated annealing algorithm.* Let $\lambda > 0$ be given. Since the Glauber dynamics mixes in time $O(n \log n)$ whenever $\lambda < 2/(\Delta - 2)$, we assume without loss of generality that $\lambda \geq 2/(\Delta - 2) \geq 1/\Delta > 1/3n$. Define a sequence $\lambda_0 < \lambda_1 < \cdots < \lambda_k$, by $\lambda_0 = 0$, $\lambda_k = \lambda$, and for $1 \leq i \leq k-1$, $\lambda_i = (1 + 1/3n)^{i-1}/3n$, where $k = \lceil \log(3n\lambda)/\log(1 + 1/3n) \rceil$.

We use the following simulated annealing algorithm. Let $Y_0$ be the empty set. For $1 \leq i \leq k$, simulate the Glauber dynamics at fugacity $\lambda_i$ for $T_i = \Omega(n \log n)$ steps, starting from initial state $Y_{i-1}$, and let $Y_i$ be the state reached after $T_i$ steps. Let the constant hidden in the $\Omega$ notation be large enough that Lemma 4.1 guarantees that Glauber dynamics mixes to within $\delta/k$ of the Gibbs distribution at fugacity $\lambda_i$, within $T$ steps, from any warm start. Output the final state $Y_k$.

The proof of correctness, presented in the next section, relies on showing that, for $1 \leq i \leq k$, the Gibbs distribution at fugacity $\lambda_{i-1}$ is a warm start to the Gibbs distribution at fugacity $\lambda_i$. Once this has been established, Lemma 4.1 will complete the proof.

PROOF OF THEOREM 1.5. Assume without loss of generality that $\lambda \geq 1/\Delta > 1/n$. Note, for $\lambda < 1/\Delta$ there is a straightforward coupling argument which proves the Glauber dynamics is close to its stationary distribution after $O(n \log n)$ steps; see [18] for a more complicated argument when $\lambda < 2/(\Delta - 2)$.

First, we prove that, for $1 \leq i \leq k$, the Gibbs distribution at fugacity $\lambda_{i-1}$ is a warm start to the Gibbs distribution at fugacity $\lambda_i$. For each $1 \leq i \leq k$, define the "partition function" $Z_i$ by

$$Z_i = \sum_{\sigma \in \Omega} \lambda_i^{|\sigma|}.$$

It is clear that the desired warm-start condition is equivalent to $Z_i \leq 2Z_{i-1}$.

We handle the case $i = 1$ separately:

$$Z_1 = \sum_{\sigma \in \Omega} \lambda_1^{|\sigma|} \leq \sum_{0 \leq i \leq n} \binom{n}{i} \lambda_1^i = \left(1 + \frac{1}{3n}\right)^n < e^{1/3} < 2 = 2Z_0.$$

For $2 \leq i \leq k$, we use the fact that $\lambda_i \leq (1 + 1/3n)\lambda_{i-1}$. From this it immediately follows that

$$Z_i = \sum_{\sigma \in \Omega} \lambda_i^{|\sigma|} \leq \sum_{\sigma \in \Omega} \lambda_{i-1}^{|\sigma|}(1 + 1/3n)^{|\sigma|} \leq Z_{i-1}(1 + 1/3n)^n < e^{1/3}Z_{i-1}.$$

This establishes the warm-start condition.

Now, for $0 \leq i \leq k$, let $\pi_i$ denote the Gibbs distribution for fugacity $\lambda_i$. We now prove by induction on $i$ that $\|Y_i - \pi_i\|_{\mathrm{TV}} \leq i\delta/k$, which will complete the proof.

The base case $i = 0$ is trivial. Let $i \geq 1$, and suppose by inductive hypothesis, $\|Y_{i-1} - \pi_{i-1}\|_{\mathrm{TV}} \leq (i-1)\delta/k$. To understand the distribution of $Y_i$, let $T = T_i$, and let us examine a $T$-step coupling of two copies of the Glauber dynamics at fugacity $\lambda_i$. Sample $A_0$ from the distribution of $Y_{i-1}$, and $B_0$ according to $\pi_{i-1}$; couple these distributions so that

$$\mathbf{Pr}(A_0 \neq B_0) = \|Y_0 - \pi_{i-1}\|_{\mathrm{TV}} \leq (i-1)\delta/k.$$

Sample $A_1, B_1, \ldots, A_T, B_T$ using a maximal coupling of the dynamics at fugacity $\lambda_i$. Since $B_0$ was a warm start, Lemma 4.1 tells us $\|B_T - \pi_i\|_{\mathrm{TV}} \leq \delta/k$. Now the triangle inequality tells us

$$\|Y_i - \pi_i\|_{\mathrm{TV}} = \|A_T - \pi_i\|_{\mathrm{TV}} \leq \mathbf{Pr}(A_T \neq B_T) + \|B_T - \pi_i\|_{\mathrm{TV}}$$
$$\leq \mathbf{Pr}(A_0 \neq B_0) + \|B_T - \pi_i\|_{\mathrm{TV}} \leq i\delta/k. \qquad \square$$

4.1. *Local uniformity.* The main result of this section is a local uniformity property of random independent sets. As before, we assume the graph is $\Delta$-regular and the fugacity, $\lambda$, is less than $e/\Delta$. For convenience, we will also assume $\lambda \geq 1/\Delta$; we do not think this condition is necessary, but it will simplify the proof of one of our technical results, Lemma 4.6. Stated roughly, the uniformity property is that every vertex has about the same number of "unblocked" neighbors, by which we mean neighbors which could be added to the independent set without violating the independence condition, and moreover that this is a fairly small fraction of $\Delta$. In Section 4.2 we will show that any pair of independent sets with this uniformity property is distance-decreasing. These two results, together with Theorem 1.3, will be the key ingredients in the proof of Lemma 4.1, given in Section 4.3.

We will use the following notation. For an independent set $X \subseteq V$ and vertex $v$, let

$$U(X, v) := \{w \in N(v) : X \cap N^*(w) = \varnothing\},$$

where

$$N^*(w) = N_v^*(w) = N(w) \setminus \{v\}.$$

Thus, $U(X, v)$ denotes the set of neighbors of $v$ that are unblocked in $X \setminus \{v\}$. For real numbers $a, b$ with $b \geq 0$, we will use the shorthand $a \pm b$ to denote the interval $[a - b, a + b]$.

LEMMA 4.2. *Let $\zeta > 0$, let $G = (V, E)$ be a $\Delta$-regular graph of girth $\geq 6$ and let $1/\Delta \leq \lambda \leq (1 - \zeta)e/\Delta$. Let $\mu$ be the solution to $\mu = \exp(-\mu\lambda\Delta)$. Let $X \subseteq V$ be a random independent set drawn from the Gibbs distribution $\pi$ at fugacity $\lambda$. Then for every $\xi > 0$,*

$$\mathbf{Pr}((\forall\, v \in V)|U(X, v)| \in (\mu \pm \xi)\Delta)$$
$$\geq 1 - 3n\exp\left(-\left(\frac{\xi\zeta}{8\lambda\Delta} - \frac{(e+1)^2}{\Delta}\right)^2 \frac{\Delta}{8}\right).$$

*In particular, for any fixed values $\zeta, \xi, \delta > 0$, there exists $C > 0$ such that $\Delta \geq C \log n$ implies*

$$\mathbf{Pr}((\forall v \in V)|U(X,v)| \in (\mu \pm \xi)\Delta) \geq 1 - \delta/n^2.$$

We conjecture that a similar uniformity property should be true for non-regular graphs of maximum degree $\Delta$; however, in this setting, $\mu$ would not be a constant, but rather a function of the vertex $v$. If so, the rest of our results would also extend easily to this context.

Our proof of Lemma 4.2 can be modified to give nontrivial bounds even when $\Delta = O(1)$. In this range, one must avoid taking union bounds over the vertex set. However, our proof of Lemma 4.1 requires such a union bound in another step, so we have focused on the case $\Delta = \Omega(\log n)$, which simplifies our results and proofs.

The first tool for our proof of Lemma 4.2 is a "bootstrapping" mechanism, to convert local recurrences for functions on $V$ into absolute bounds.

We will use the following notational conventions.

DEFINITION 4.3. Let $\mathscr{I}$ denote the set of all real intervals, $[a,b]$. For any set $S \subseteq \mathbb{R}$, let $S \pm \xi$ denote the set $S + [-\xi, \xi] = \{x + y \mid x \in S, |y| < \xi\}$.

LEMMA 4.4. *Let $G = (V, E)$ be a graph, and let $f: V \to [0,1]$ be a random function (from any distribution over $[0,1]^V$). Let $\theta, \varphi > 0$, let $g: \mathbb{R} \to \mathbb{R}$ and define $h: \mathscr{I} \to \mathscr{I}$ by $h(I) = g(I) \pm \theta$. Suppose that for every vertex $v$,*

$$\mathbf{Pr}\left(f(v) \notin \frac{1}{\Delta}\sum_{w \sim v} g(f(w)) \pm \theta\right) \leq \varphi.$$

*Then*

$$\mathbf{Pr}((\forall v \in V, \forall t \geq 1)\, f(v) \in h^t([0,1])) \geq 1 - n\varphi.$$

PROOF. A union bound over $v \in V$ shows that

(4) $$\mathbf{Pr}\left((\forall v \in V)\, f(v) \in \frac{1}{\Delta}\sum_{w \in N(v)} g(f(w)) \pm \theta\right) \geq 1 - n\varphi.$$

Suppose this event holds, but that there exists some vertex $z$ and positive integer $t$ such that $f(z) \notin h^t([0,1])$. Without loss of generality, choose $z$ so that $t$ is as small as possible. Then in particular, for all $w \sim z$, $f(w) \in h^{t-1}([0,1])$. But in this case, specializing (4) to vertex $z$ implies

$$f(z) \in \frac{1}{\Delta}\sum_{w \sim z} h(w) \subset h^t([0,1]),$$

a contradiction. □

The next step in proving Lemma 4.2 is the following easy observation, which we will use again in Section 4.3.

OBSERVATION 4.5.　*Let $\zeta \in (0,1)$, let $C = (1-\zeta)e$ and let $\mu = \exp(-C\mu)$. Then*

$$\mu < (1 - \zeta/2)/C.$$

PROOF.　Let $f(x) = \exp(-Cx)$ and let $y = (1 - \zeta/2)/C$. Then

$$\begin{aligned}\frac{f(y)}{y} &= \frac{\exp(\zeta/2)}{ey} \\ &= \frac{(1-\zeta)\exp(\zeta/2)}{1 - \zeta/2} \\ &< (1 - \zeta/2)\exp(\zeta/2) < 1,\end{aligned}$$

where the last inequality is since $1 - x < \exp(-x)$ for all $x \neq 0$. Thus $f(y) < y$, which implies $y > \mu$, since $f$ is a decreasing function with fixed point $\mu$.　□

Our next result is a strong sort of convergence for iterated applications of the mapping $x \mapsto \exp(-Cx)$, where $1 \le C < e$. It says that, even permitting a small adversarial perturbation after every step, every trajectory under the iterated mapping quickly converges to a small interval around the unique fixed point $\mu = \exp(-C\mu)$.

LEMMA 4.6.　*Let $1 \le C \le (1-\zeta)e$ and $\xi > 0$. Let $g(x) = \exp(-Cx)$, let $\mu$ be the unique fixed point of $g$ and set $\theta = \xi\zeta/8C$. Let $h(I) = g(I) \pm \theta$, and let $t = \lceil \frac{4}{\zeta}\ln(1 + \frac{1}{\xi}) \rceil$. Then*

$$h^t([0,1]) \subseteq \mu \pm \xi.$$

REMARK 4.7.　The assumption $C \ge 1$ is only for convenience in our proof; since rapid mixing is already known when $C \le 2$, there was no particular reason to handle the case $0 < C < 1$. However, the assumption $C \le e$ is necessary. Indeed, for $C > e$, the unique fixed point for $x \mapsto \exp(-Cx)$ is unstable, and from any starting point $x_0 \neq \mu$, the sequence $g^i(x_0)$ approaches a fixed cycle of period 2 (also unique).

PROOF OF LEMMA 4.6.　Without loss of generality, $\zeta, \xi \le 1$; otherwise there is nothing to prove. Define $I : [0,1] \to \mathscr{I}$ by

$$I(x) = \left[\mu - \frac{x}{C}, \mu - \frac{\ln(1-x)}{C}\right],$$

where we adopt the convention that $I(1) = [\mu - 1/C, \infty)$. By the definition of $h$, since $g$ is decreasing, and since $\mu$ is the fixed point of $g$, we have

$$h(I(x)) = [\mu(1-x) - \theta, \mu e^x + \theta].$$

Next we will prove that, whenever $x \geq \xi/2$,

(5) $$h(I(x)) \subset I(x(1 - \zeta/4)).$$

This is equivalent to checking that

(6) $$\mu x + \theta \leq \frac{x(1 - \zeta/4)}{C}$$

and

(7) $$\mu(e^x - 1) + \theta \leq -\frac{\ln(1 - x(1 - \zeta/4))}{C}.$$

Since $x \geq \xi/2$, we have $\theta = \zeta\xi/8C \leq \zeta x/4C$. By Observation 4.5, we also know that $\mu \leq (1 - \zeta/2)/C$. This allows us to deduce the first desired inequality, thus

$$\mu x + \theta \leq \frac{(1 - \zeta/2)x}{C} + \frac{\zeta x/4}{C}$$
$$= \frac{x(1 - \zeta/4)}{C}.$$

To deduce the second inequality, we make the same substitutions for $\mu$ and $\theta$; then we show that the inequality holds termwise for the Taylor expansions around the origin:

$$\theta + \mu(e^x - 1) \leq \frac{\zeta x}{4C} + \frac{(1 - \zeta/2)}{C}(e^x - 1)$$
$$= \frac{\zeta x}{4C} + \frac{(1 - \zeta/2)}{C} \sum_{j \geq 1} \frac{x^j}{j!}$$
$$= \frac{(1 - \zeta/4)x}{C} + \frac{(1 - \zeta/2)}{C} \sum_{j \geq 2} \frac{x^j}{j!}$$
$$\leq \frac{(1 - \zeta/4)x}{C} + \frac{1}{C} \sum_{j \geq 2} \frac{(1 - \zeta/4)^j x^j}{j}$$
$$= \frac{1}{C} \sum_{j \geq 1} \frac{(1 - \zeta/4)^j x^j}{j}$$
$$= \frac{-\ln(1 - x(1 - \zeta/4))}{C}.$$

Thus we have established $h(I(x)) \subset I(x(1-\zeta/4))$ whenever $x \geq \xi/2$. It follows by induction that for any positive integer $k$, as long as $x(1-\zeta/4)^{k-1} \geq \xi/2$,

$$h^k(I(x)) \subset I(x(1-\zeta/4)^k).$$

In particular, let $k$ be the least positive integer such that $(1-\zeta/4)^k \leq C\xi/(1+C\xi)$. Then $h^k(I(1)) \subset I(C\xi/(1+C\xi))$. By another application of Observation 4.5, we may deduce $[0,1] \subseteq I(1) = [\mu - 1/C, \infty)$. Also, elementary algebra implies $I(C\xi/(1+C\xi)) \subset \mu \pm \xi$. It follows by monotonicity of $h$ that $h^k([0,1]) \subset h^k(I(1)) \subset I(C\xi/(1+C\xi)) \subset \mu \pm \xi$.

Solving for $k$, we find

$$k = \left\lceil \frac{\ln(C\xi/(1+C\xi))}{\ln(1-\zeta/4)} \right\rceil$$

$$\leq \left\lceil \frac{4}{\zeta} \ln((1+C\xi)/C\xi) \right\rceil$$

$$\leq \left\lceil \frac{4}{\zeta} \ln(1+1/\xi) \right\rceil,$$

which completes the proof. □

Now we are equipped to prove Lemma 4.2.

PROOF OF LEMMA 4.2. Fix a vertex $v \in V$. For each neighbor $w$ of $v$, let $Y_w$ denote the indicator variable for the event that $X \cap N^*(w) = \varnothing$. Then

$$|U(X,v)| = \sum_{w \sim v} Y_w.$$

Now, for each neighbor $w$, let us compute the conditional expectation of $Y_w$ given $X \setminus N^*(w)$, that is, given $X \cap (V \setminus N^*(w))$. By definition, $Y_w$ equals 1 iff no neighbor of $w$ is in $X$, except possibly $v$. Hence if $w \in X$, we have

$$\mathbf{E}(Y_w | X \setminus N^*(w)) = 1.$$

If $w \notin X$, then

$$\mathbf{E}(Y_w | X \setminus N^*(w)) = \prod_{z \in U(X,w) \setminus \{v\}} 1/(1+\lambda)$$

$$= (1+\lambda)^{-|U(X,w) \setminus \{v\}|}$$

$$\in (e^{-\lambda}, e^{-\lambda+\lambda^2})^{|U(X,w) \setminus \{v\}|}$$

$$= (1, e^{\lambda^2 |U(X,w) \setminus \{v\}|}) e^{-\lambda |U(X,w) \setminus \{v\}|}$$

$$\subset (1, e^{\lambda^2(\Delta-1)}) e^{-\lambda |U(X,w) \setminus \{v\}|}$$

$$\subset (1, e^{e\lambda})e^{-\lambda|U(X,w)\setminus\{v\}|}$$
$$\subset (1, e^{(e+1)\lambda})e^{-\lambda|U(X,w)|}$$
$$\subset e^{-\lambda|U(X,w)|} \pm (e^{(e+1)\lambda} - 1)$$

assuming $e^{(e+1)\lambda} \leq e^{(e+1)e/\Delta} \leq 1 + (e^2 + e + 1)/\Delta$, which is true whenever $\Delta \geq 100$. Since the desired upper bound on probability is trivial for $\Delta \leq 750$, the above assumption may be made without loss of generality.

Let $S_2(v)$ denote the set of vertices at distance 2 from $v$, that is,

$$S_2(v) = \bigcup_{w \in N(v)} N^*(w).$$

Applying linearity of expectation and then averaging, we have

$$\mathbf{E}(|U(X,v)||X \setminus S_2(v))$$
$$= \sum_{i=1}^{\Delta} \mathbf{E}(Y_w | X \setminus S_2(v))$$
$$\in |X \cap N(v)| + \left(\sum_{w \in N(v)\setminus X} \exp(-\lambda|U(X,w)|)\right) \pm (e^2 + e + 1).$$

Since the girth is at least 6, there are no edges between vertices in $S_2(v)$. Hence, conditioned on $X \setminus S_2(v)$, the random variables $Y_w$ are fully independent, and take values in $[0,1]$. It follows by Chernoff's bound that, for all $\psi > 0$,

$$\mathbf{Pr}\bigg(|U(X,v)| \notin |X \cap N(v)|$$

(8)
$$+ \sum_{w \in N(v)\setminus X} \exp(-\lambda|U(X,w)|) \pm (e^2 + e + 1 + \psi\Delta/2)\bigg)$$
$$\leq 2\exp(-\psi^2\Delta/8).$$

If $v \in X$, we have $|X \cap N(v)| = 0$. When $v \notin X$, since $\lambda < e/\Delta$,

$$\mathbf{E}(|X \cap N(v)||X \setminus N(v)) \leq \Delta \frac{\lambda}{1+\lambda} < e.$$

Apply Chernoff's bound again, this time to the sum of indicator variables for the events $\{w \in X\}$, where $w$ is a neighbor of $v$. These events are conditionally independent, given $X \setminus N(v)$, and so, since the uniform upper bound of $e$ applies to the conditional expectation, Chernoff's bound implies, for all $\psi > 0$,

(9) $$\mathbf{Pr}(|X \cap N(v)| > e + \psi\Delta/2) \leq \exp(-\psi^2\Delta/8).$$

Combining (8) and (9), we have, for all $\psi > 0$,

$$\mathbf{Pr}\left(|U(X,v)| \notin \sum_{w \sim v} \exp(-\lambda|U(X,w)|) \pm ((e+1)^2 + \psi\Delta)\right) \leq 3\exp(-\psi^2\Delta/8).$$

Applying Lemmas 4.4 and 4.6 with parameters $f(v) = |U(X,v)|/\Delta$, $C = \lambda\Delta$, $g(x) = \exp(-Cx)$, $\theta = \xi\zeta/8C$, $h = g \pm \theta$ and $\varphi = 3\exp(-(\theta - (e+1)^2/\Delta)^2\Delta/8)$, there must exist some $t \geq 1$ such that

$$\mathbf{Pr}((\forall v \in V)|U(X,v)| \notin (\mu \pm \xi)\Delta) \leq \mathbf{Pr}((\forall v \in V)\,|U(X,v)| \notin h^t([0,1]))$$
$$\leq n\varphi. \qquad \square$$

4.2. *Convergence of the coupling.* We now show that pairs $(X, Y)$ are distance-decreasing, so long as no vertex of $G$ has too many unblocked neighbors with respect to either set.

LEMMA 4.8. *Let $G$ be a graph, and consider the Glauber dynamics for independent sets, with fugacity $\lambda$. Let $X$ and $Y$ be independent sets and suppose there exists $\zeta > 0$ such that for every vertex $v$,*

$$|U(X,v)|, |U(Y,v)| \leq (1-\zeta)\frac{1+\lambda}{\lambda}.$$

*Then the pair $(X, Y)$ is $\zeta/n$ distance-decreasing.*

PROOF. Let $X', Y'$ be the new sets obtained after doing one step of Glauber dynamics starting from $(X, Y)$, using a maximal coupling. More explicitly, select a uniformly random vertex $v^*$ for update. If the neighborhood of $v^*$ is disjoint from both $X$ and $Y$, then with probability $\lambda/(1+\lambda)$, set $X' = X \cup \{v^*\}$ and $Y' = Y \cup \{v^*\}$, and otherwise set $X' = X$ and $Y' = Y$. If the neighborhood of $v^*$ is disjoint from exactly one of $X$ or $Y$, add $v^*$ to the corresponding new set $X'$ or $Y'$ with probability $\lambda/(1+\lambda)$. Otherwise, make no change.

Let $D = \{v : X(v) \neq Y(v)\}$, and $D' = \{v : X'(v) \neq Y'(v)\}$. Then, letting $\rho$ denote Hamming distance, we have

$$\mathbf{E}(\rho(X',Y')) - \rho(X,Y)$$
$$= \mathbf{E}(|D'|) - |D|$$
$$= -\mathbf{Pr}(v^* \in D) + \mathbf{Pr}(v^* \in D')$$
$$\leq -\frac{\rho(X,Y)}{n} + \sum_{w \in D}\sum_{v \in N(w)} \mathbf{Pr}(v^* = v \text{ and } X'(v) \neq Y'(v))$$
$$= -\frac{\rho(X,Y)}{n} + \sum_{w \in X \setminus Y} |U(Y,w)|\frac{\lambda}{n(1+\lambda)} + \sum_{w \in Y \setminus X} |U(X,w)|\frac{\lambda}{n(1+\lambda)}$$

$$\leq \frac{-\zeta \rho(X,Y)}{n},$$

where the last inequality follows by the hypothesis of the lemma. □

4.3. *Proof of Lemma* 4.1. We are now ready to prove Lemma 4.1, establishing rapid mixing of the Glauber dynamics from a warm start.

PROOF OF LEMMA 4.1. Let $\mu = \exp(-\mu\lambda\Delta)$. Let $S$ denote the set of independent sets $X$ such that for every vertex $v$,

$$|U(X,v)| \leq (\mu + \zeta/8)\Delta.$$

Applying Lemma 4.2 with $\xi = \zeta/8$, then simplifying using our assumptions $\lambda\Delta < e$ and $\Delta \geq 320000\ln(144n^3/\zeta\delta)/\zeta^4$, we obtain

$$\pi(S) \geq 1 - 3n\exp\left(-\left(\frac{\zeta^2}{64\lambda\Delta} - \frac{(e+1)^2}{\Delta}\right)^2 \frac{\Delta}{8}\right)$$

$$\geq 1 - 3n\exp\left(\frac{-\zeta^4\Delta}{320000}\right)$$

$$\geq 1 - \frac{\zeta\delta}{48n^2}.$$

Let $X \in S$, and let $v$ be any vertex. Applying Observation 4.5 to the fixed point $\mu$ (with $C = \lambda\Delta$), and recalling our hypothesis that $\lambda < e/\Delta$,

$$|U(X,v)| \leq (\mu + \zeta/8)\Delta$$

$$\leq \left(\frac{1-\zeta/2}{\lambda\Delta} + \zeta/8\right)\Delta$$

$$= \frac{1-\zeta/2 + \zeta\lambda\Delta/8}{\lambda}$$

$$< \frac{1-\zeta/8}{\lambda}.$$

Hence, by Lemma 4.8, every pair $(X,Y) \in S \times S$ is $\zeta/8n$ distance-decreasing.

The desired result now follows by applying the coupling with stationarity Theorem 1.3 to the set $S$. □

**Acknowledgments.** We would like to thank the anonymous referees for their careful proofreading and many helpful suggestions.

## REFERENCES

[1] ALDOUS, D. J. (1983). Random walks on finite groups and rapidly mixing Markov chains. *Séminaire de Probabilities XVII. Lecture Notes in Math.* **986** 243–297. Springer, Berlin. MR0770418

[2] VAN DEN BERG, J. and STEIF, J. E. (1994). Percolation and the hard-core lattice gas model. *Stochastic Process. Appl.* **49** 179–197. MR1260188

[3] BUBLEY, R. and DYER, M. E. (1997). Path coupling: A technique for proving rapid mixing in Markov chains. In *Proceedings of the 38th Annual IEEE Symposium on Foundations of Computer Science* 223–231. IEEE Computer Society Press, Los Alamitos, CA.

[4] DOEBLIN, W. (1938). Exposé de la théorie des chaînes simples constantes de Markov à un nombre fini d'états. *Revue Mathématique de l'Union Interbalkanique* **2** 77–105.

[5] DYER, M. and FRIEZE, A. (2003). Randomly colouring graphs with lower bounds on girth and maximum degree. *Random Structures Algorithms* **23** 167–179. MR1995689

[6] DYER, M., FRIEZE, A., HAYES, T. and VIGODA, E. (2004). Randomly coloring constant degree graphs. In *Proceedings of the 45th Annual IEEE Symposium on Foundations of Computer Science* 582–589. IEEE Computer Society Press, Los Alamitos, CA.

[7] DYER, M., SINCLAIR, A., VIGODA, E. and WEITZ, D. (2004). Mixing in time and space for lattice spin systems: A combinatorial view. *Random Structures Algorithms* **24** 461–479. MR2060631

[8] GOLDBERG, L. A., MARTIN, R. and PATERSON, M. (2004). Strong spatial mixing for graphs with fewer colours. In *Proceedings of the 45th Annual IEEE Symposium on Foundations of Computer Science* 562–571. IEEE Computer Society Press, Los Alamitos, CA.

[9] HAYES, T. P. (2003). Randomly coloring graphs of girth at least five. In *Proceedings of the 35th Annual ACM Symposium on Theory of Computing* 269–278. ACM, New York. MR2121038

[10] HAYES, T. P. and VIGODA, E. (2003). A non-Markovian coupling for randomly sampling colorings. In *Proceedings of the 44th Annual IEEE Symposium on Foundations of Computer Science* 618–627. IEEE Computer Society Press, Los Alamitos, CA.

[11] JERRUM, M. R. (1995). A very simple algorithm for estimating the number of $k$-colourings of a low-degree graph. *Random Structures Algorithms* **7** 157–165. MR1369061

[12] JERRUM, M. R., SINCLAIR, A. and VIGODA, E. (2004). A polynomial-time approximation algorithm for the permanent of a matrix with non-negative entries. *J. Assoc. Comput. Machinery* **51** 671–697. MR2147852

[13] KANNAN, R., LOVÁSZ, L. and SIMONOVITS, M. (1997). Random walks and an $O^*(n^5)$ volume algorithm for convex bodies. *Random Structures Algorithms* **11** 1–50. MR1608200

[14] LINDVALL, T. (2002). *Lectures on the Coupling Method*. Dover Publications Inc., New York. MR1924231

[15] MARTINELLI, F. (1999). Lectures on Glauber dynamics for discrete spin models. *Lectures on Probability Theory and Statistics. Lecture Notes in Math.* **1717** 93–191. Springer, Berlin. MR1746301

[16] MOLLOY, M. (2004). The Glauber dynamics on colorings of a graph with high girth and maximum degree. *SIAM J. Comput.* **33** 712–737. MR2066651

[17] VIGODA, E. (2000). Improved bounds for sampling colorings. *J. Math. Phys.* **41** 1555–1569. MR1757969

[18] VIGODA, E. (2001). A note on the Glauber dynamics for sampling independent sets. *Electron. J. Combin.* **8** 1–8. MR1814515

Computer Science Division  
University of California at Berkeley  
Berkeley, California 94720  
USA  
E-mail: hayest@cs.berkeley.edu

College of Computing  
Georgia Institute of Technology  
Atlanta, Georgia 30332  
USA  
E-mail: vigoda@cc.gatech.edu